\documentclass[twoside,leqno]{article}

\usepackage[letterpaper]{geometry}

\usepackage{amsmath}
\usepackage{amssymb}
\usepackage{amsthm}
\usepackage[T1]{fontenc}
\usepackage{amsfonts}
\usepackage{graphicx}
\usepackage{epstopdf}
\usepackage{enumitem}
\usepackage{algorithmic}
\usepackage{ifpdf}
\usepackage{url}

\ifpdf
  \DeclareGraphicsExtensions{.eps,.pdf,.png,.jpg}
\else
  \DeclareGraphicsExtensions{.eps}
\fi

\newcommand{\RR}{\mathbb{R}}
\renewcommand{\d}{\mathrm{d}}

\newcommand{\Hh}{\mathcal{H}}
\newcommand{\Pp}{\mathcal{P}}
\newcommand{\Alpha}{{\boldsymbol{\alpha}}}
\newcommand{\email}[1]{\texttt{#1}}

\newtheorem{theorem}{Theorem}
\newtheorem{proposition}{Proposition}

\theoremstyle{definition}

\theoremstyle{remark}

\usepackage{amsopn}

\begin{document}

\title{\Large Optimal and Diffusion Transports in Machine Learning}
\author{Gabriel Peyr\'e\thanks{CNRS and ENS PSL Universit\'e 
(\email{gabriel.peyre@ens.fr}, \texttt{https://www.gpeyre.com}).}}
\date{}

\maketitle

\begin{abstract} 
Several problems in machine learning are naturally expressed as the design and analysis of time-evolving probability distributions. This includes sampling via diffusion methods, optimizing the weights of neural networks, and analyzing the evolution of token distributions across layers of large language models. While the target applications differ (samples, weights, tokens), their mathematical descriptions share a common structure. A key idea is to switch from the Eulerian representation of densities to their Lagrangian counterpart through vector fields that advect particles. This dual view introduces challenges, notably the non-uniqueness of Lagrangian vector fields, but also opportunities to craft density evolutions and flows with favorable properties in terms of regularity, stability, and computational tractability. This survey presents an overview of these methods, with emphasis on two complementary approaches: diffusion methods, which rely on stochastic interpolation processes and underpin modern generative AI, and optimal transport, which defines interpolation by minimizing displacement cost. We illustrate how both approaches appear in applications ranging from sampling and neural network optimization to modeling the dynamics of transformers for large language models.
\end{abstract}

\section{Introduction.}

The goal of this survey is to explain the connection between optimal transport and evolutions over the space of probability measures, with emphasis on partial differential equations and generative models based on diffusion. The exposition is informal, favoring intuition over detailed proofs. We work over $\mathbb{R}^d$ and highlight links with machine learning applications: the training dynamics of multilayer perceptrons, the modeling of transformers where tokens are transported, and flow matching approaches for generative models.
This paper is organized around four core examples that illustrate how dynamical systems of probability measures connect with modern machine learning:  
\begin{itemize}
    \item \textbf{Diffusion models for generative AI} (\S\ref{sec:genmodels}). Generative models construct flows $(\alpha_t)_t$ between simple references and complex data distributions. Diffusion and flow matching models achieve state-of-the-art results, but raise fundamental questions on sample complexity, discretization errors, and geometric interpretation.  

    \item \textbf{Optimal transport and ML} (\S\ref{sec:ot}). Rooted in Monge and Kantorovich, optimal transport provides a geometric framework for comparing and transforming measures. It has become a central tool in statistics and machine learning, with scalable entropic algorithms and dynamic formulations (Benamou--Brenier) enabling practical use in generative modeling and beyond.  

    \item \textbf{Optimization over measures} (\S\ref{sec:wass-flow} and \S\ref{sec:mlp}). Wasserstein gradient flows extend OT to optimization, interpreting nonlinear PDEs as steepest descents in Wasserstein space. This theory informs both PDE analysis (porous medium, crowd motion) and machine learning (training shallow neural networks as mean-field limits with convergence guarantees).  

    \item \textbf{Very deep transformers} (\S\ref{sec:transformers}). Transformer architectures process data as sets of tokens coupled through self-attention. Their depth can be modeled by continuous-time flows or PDEs over token distributions, linking attention dynamics to Vlasov- or aggregation-type equations. This measure-theoretic view sheds light on expressivity, clustering, and universality.  
\end{itemize}

Taken together, these examples reveal a common structure: learning systems can often be described as evolutions of probability measures under vector fields, sometimes with a variational or gradient-flow interpretation, sometimes not. What is at stake is a unified mathematical understanding of how such flows encode computation, optimization, and generalization in machine learning.

\section{Evolutions over the Space of Measures}

All the problems outlined in the introduction can be formulated by studying the evolution $t \mapsto \alpha_t \in \mathcal{P}(\mathbb{R}^d)$, where $\alpha_t$ is a probability measure on $\mathbb{R}^d$ with unit mass. Such evolution can be described in a ``Lagrangian'' way as the advection of particles along a (time-dependent) vector field $v_t(x)$ in $\mathbb{R}^d$. At the particle level, this advection is governed by 
\begin{equation}
    \frac{\mathrm{d}x(t)}{\mathrm{d}t} = v_t(x(t)), \label{eq:lagrangian-advection}
\end{equation}
such that $x(0)$ is mapped to $x(t)$ by a ``transport'' mapping $T_t : x(0) \mapsto x(t)$. 
The fact that $\alpha_t$ is the law of the advected particles implies $\alpha_t = (T_t)_\sharp \alpha_0$. 
Here $T_\sharp$ is the ``push-forward'' operator between measures~\cite{santambrogio2015optimal}, which for discrete measures acts by transporting the support of the Dirac masses, i.e. $T_\sharp\!\left( \tfrac{1}{n} \sum_i \delta_{x_i} \right) = \tfrac{1}{n} \sum_i \delta_{T(x_i)}$, and this definition extends to arbitrary measures (for instance with densities). This means that for the flow we consider, if $\alpha_0 = \tfrac{1}{n} \sum_i \delta_{x_i(0)}$ represents a discrete point cloud, then $\alpha_t = \tfrac{1}{n} \sum_{i=1}^n \delta_{x_i(t)}$, where each $x_i(t)$ solves \eqref{eq:lagrangian-advection}.

Another way to describe this evolution is to track the measure $\alpha_t$ itself. For instance, one may consider its density $\rho_t = \tfrac{\d \alpha_t}{\d x}$ with respect to the Lebesgue measure, provided it exists. In this case, $\rho_t(x)$ represents the local concentration of particles around $x$. In this ``Eulerian'' representation, the ODE for particle trajectories becomes the PDE~\cite{ambrosio2008gradient}
\begin{equation} \label{eq:eulerian-advection}
	\frac{\partial \alpha_t}{\partial t} + \mathrm{div}(v_t \alpha_t) = 0,
\end{equation} 
where $\mathrm{div}(v_t \alpha_t)$ is understood in the weak sense, or in terms of the density $\rho_t$ whenever it exists. This PDE is commonly called the advection equation, the continuity equation, or Liouville's equation when defined on a phase space. Finally, if $\alpha_0 = \tfrac{1}{n} \sum_i \delta_{x_i(0)}$ is a discrete particle distribution, then $\alpha_t = \tfrac{1}{n} \sum_i \delta_{x_i(t)} = (T_t)_\sharp \alpha_0$ is again a sum of moving particles with velocities $\dot x_i(t) = v_t(x_i(t))$, and it satisfies the PDE only in the weak sense.

\subsection{From measure evolutions to vector fields.}

A crucial point, central to both the practical and theoretical aspects of this review, is that for a given evolution $\alpha_t$, there are infinitely many possible choices of vector fields $v_t$ satisfying~\eqref{eq:eulerian-advection}. This is because adding to $v_t$ any time-dependent field $w_t$ from the linear space $\Hh_{\alpha_t} := \{ w : \mathrm{div}(\alpha_t w) = 0 \}$ does not affect the evolution of the measure. In other words, $\Hh_\alpha$ consists of vector fields that leave $\alpha$ invariant. For example, if $\alpha_0$ is an isotropic measure such as a Gaussian $\mathcal{N}(0,\mathrm{Id}_d)$, then orthogonal rotations $T_t$ (and in fact a much richer class of transformations) leave $(T_t)_\sharp \alpha_0 =  \alpha_0$ unchanged. More generally, $\Hh_\alpha$ is non-trivial since it corresponds to the kernel of a weighted divergence operator. In the Gaussian case, $\Hh_\alpha$ includes vector fields generated by antisymmetric matrices, which induce rotations. These ambiguities are both a challenge—since they require some form of regularization to make the problem well-posed—and an opportunity, as they allow the design of vector fields with desirable computational and statistical properties.

\paragraph{Dacorogna and Moser inversion.}

Reconstructing particle evolution, i.e., recovering a vector field $v_t$, from an observed measure evolution $t\mapsto \alpha_t$ (sometimes called ``trajectory inference'') is thus an ill-posed inverse problem.
A simple choice is to impose that $\alpha_t v_t$ is a gradient field, thus leading to the inversion of a Laplacian (which is possible assuming boundary conditions, for instance, vanishing at infinity):
\begin{equation}\label{eq:dacorogna-moser}
    v_t = - \frac{1}{\alpha_t} \nabla \Delta^{-1}(\partial_t \alpha_t).
\end{equation}
This construction was initially proposed by Dacorogna and Moser~\cite{DacorognaMoser1990}. A difficulty with this choice is that it is not well defined when $\alpha_t$ vanishes. Moreover, the resulting velocity is not in general a gradient field, a property that will be desirable in several cases below. In the following, we consider techniques where $v_t$ is a gradient by design, is well defined, and can be computed more efficiently without explicitly inverting a Laplacian.  

\paragraph{Least square inversion and gradient structure.}
\label{sec:least-square-field}

A fruitful approach is to view the reconstruction of $v_t$ from $\alpha_t$ as an inverse problem and to introduce regularization. A standard choice is to use a least-squares formulation, that is, to impose a quadratic prior on the recovered vector field. Since we are advecting probability distributions, the quadratic norm over vector fields should be weighted by the measure itself, leading to the following quadratic problem under linear constraints:
\begin{equation}\label{eq:least-square-field}
	\min_v \frac{1}{2}\int_0^1 \int_{\mathbb{R}^d} \|v_t(x)\|^2 \, \mathrm{d}\alpha_t(x) \, \mathrm{d}t \quad \text{subject to} \quad \mathrm{div}(\alpha_t v_t) + \partial_t \alpha_t = 0.
\end{equation}
The idea of introducing a measure-dependent quadratic norm over the space of vector fields can be interpreted, at an abstract level, as endowing the space of probability measures with a Riemannian structure. This was first proposed by Lafferty~\cite{lafferty1988density}, and was later developed rigorously by Otto~\cite{otto2000generalization}, who studied the geometry of probability distributions in this framework. This viewpoint is central to the connection with optimal transport discussed in Section~\ref{sec:ot}.
To analyze the structure of the solution to~\eqref{eq:least-square-field}, a standard approach is to introduce a Lagrange multiplier $\psi_t(x)$, which acts as the dual variable enforcing the continuity equation. The corresponding first-order optimality conditions then establish a direct relation between $v_t$ and the potential $\psi_t$.
\begin{equation}
	v_t(x) = \nabla \psi_t(x).
\end{equation}
Hence, the optimal vector field is necessarily a gradient. This gradient structure, which will reappear in several places in the following sections (for diffusion models, optimal transport maps, and Wasserstein gradient flows), is precisely equivalent to imposing a least-effort constraint on the estimated velocity.
Substituting this back, one obtains an equation linking \(\psi_t\) and \(\partial_t \alpha_t\). This results in the explicit expression:
\begin{equation}\label{eq:least-square-field-explicit}
    v_t = -\nabla \Delta_\alpha^{-1}(\partial_t \alpha_t), \quad \text{with} \quad \Delta_\alpha(\varphi) := \mathrm{div}(\alpha_t \nabla \varphi).
\end{equation}
In general, performing this inversion is computationally demanding, but under specific choices for \(\alpha_t\), simpler formulas may be derived as we will exemplify in the following sections.

\section{Generative Models via Flow Matching}
\label{sec:genmodels}

A central challenge in generative modeling is to sample from an unknown distribution that is only accessible through a finite dataset. Modern generative approaches in AI~\cite{rezende2015variational,ho2020denoising} address this by constructing a transport map $T$ between a tractable reference distribution $\alpha_0$ (typically a standard Gaussian) and a complex target distribution $\alpha_1 = T_\sharp \alpha_0$. This enables sampling from $\alpha_1$ by first drawing $X_0 \sim \alpha_0$ and then applying the transport map (often parameterized as a neural network) so that $X_1 := T(X_0) \sim \alpha_1$. While the existence of such a map is guaranteed for any $\alpha_1$~\cite{santambrogio2015optimal}, finding an explicit and computationally efficient construction of $T$ is highly non-trivial. A recent line of work, which we present in this section, tackles this by introducing a time-dependent interpolation $(\alpha_t)_{t \in [0,1]}$ between $\alpha_0$ and $\alpha_1$, and then learning a vector field $v_t$ that satisfies the continuity equation~\eqref{eq:eulerian-advection}. For certain choices of interpolation, an admissible $v_t$ can be computed without explicitly inverting a Laplacian (as in~\eqref{eq:least-square-field-explicit} and \eqref{eq:dacorogna-moser}), but instead through a conditional expectation. This expectation can itself be estimated by solving an unconstrained least-squares problem, making the approach practical when only finite samples from $\alpha_0$ and $\alpha_1$ are available.
This general procedure has become the backbone of state-of-the-art generative models for images, audio, and video~\cite{ho2020denoising,ho2022video}.

In general, the vector field $v_t$ cannot be computed in closed form. A computationally tractable route, when Gaussian smoothing is used, is score matching~\cite{hyvarinen2005estimation}, which estimates the score of the smoothed density and hence the velocity through identities such as the one below. Building on this idea, Vincent~\cite{vincent2011connection} showed that the score is related to the optimal Gaussian denoiser under a mean-squared error criterion. This insight paved the way for the class of models now known as diffusion models~\cite{song2020score,ho2020denoising}, where sampling is achieved by integrating a time-dependent (possibly stochastic) differential equation initialized from $\alpha_0$. More recently, these methods have been extended to arbitrary interpolations between measures, giving rise to the framework of flow matching~\cite{lipman2022flow,albergo2023stochastic}, which is the focus of this section.  
Unlike optimal transport (discussed in Section~\ref{sec:ot}), the interpolations $\alpha_t$ generated by diffusion or flow matching are not defined by a displacement $T$ that is a gradient field, and their geometric characterization remains largely unresolved~\cite{lavenant2022flow}. 

\subsection{Stochastic interpolant}
\label{sec:stoch-int}

The flow-matching approach~\cite{lipman2022flow,albergo2023stochastic} constructs an explicit interpolation $(\alpha_t)_{t \in [0,1]}$ between two given distributions $(\alpha_0,\alpha_1)$, where $\alpha_0$ is typically a simple latent distribution (e.g., a Gaussian) and $\alpha_1$ is the data distribution. In its simplest form, the interpolation is defined by linearly interpolating independent random vectors 
\[
	\alpha_t := \mathrm{Law}((1-t)X_0 + tX_1), \quad\text{where}\quad 
	(X_0,X_1) \sim \alpha_0 \otimes \alpha_1.
\]
In terms of measures, if we define $S_t(x_0,x_1) := (1-t)x_0 + t x_1 \in \mathbb{R}^d$, then  
$\alpha_t := (S_t)_\sharp (\alpha_0 \otimes \alpha_1)$.  
An alternative view is that interpolating independent random vectors corresponds to convolving their rescaled laws:
\[
    \alpha_t = \mathrm{Law}((1-t)X_0) * \mathrm{Law}(tX_1),
\]
where $*$ denotes convolution. For example, when $\alpha_0$ is Gaussian, $\alpha_t$ is a blurred version of $\alpha_1$, obtained by Gaussian convolution, which is equivalent to a heat diffusion after a time reparameterization. This connection is at the heart of interpreting diffusion models, as discussed in the following section.  
More generally, one can define richer interpolations by considering couplings other than the independent product $\alpha_0 \otimes \alpha_1$, or by using nonlinear interpolation maps $S_t$.

\subsection{Flow matching field}

The central result, which underlies the practical success of this approach, is the following proposition: it provides an explicit construction of the vector field $v_t$ through a conditional expectation.

\begin{proposition}[Conditional expectation field, \cite{lipman2022flow,albergo2023stochastic}]\label{prop:fm-ce}
Let $(X_0,X_1)\sim \alpha_0 \otimes \alpha_1$ and define, for $t\in(0,1)$,
\begin{equation}\label{eq:fm-linear-ce}
  v_t(z)\;=\;\mathbb{E}\!\left[\,X_1-X_0\;\middle|\; (1-t)X_0+tX_1=z\,\right].
\end{equation}
Then the pair $(\alpha_t,v_t)$ satisfies the continuity equation~\eqref{eq:eulerian-advection} in the weak sense.
\end{proposition}

Formula~\eqref{eq:fm-linear-ce} is not directly useful in large-scale applications, since the conditional expectation is generally intractable. However, conditional expectations are precisely the optimal regressors under a mean-squared error criterion. This observation shows that $v_t$ can equivalently be computed by solving the least-squares problem
\begin{equation}\label{eq:fm-linear-obj}
  \min_{(v_t)_t}\;
  \int_0^1 \int_{\mathbb{R}^d\times\mathbb{R}^d}
  \|\,v_t\bigl((1-t)x_0+t x_1\bigr)-(x_1-x_0) \|^2\,
  \mathrm{d}\alpha_0(x_0) \mathrm{d}\alpha_1(x_1) \mathrm{d} t.
\end{equation}
This formulation has several key advantages. Unlike~\eqref{eq:least-square-field}, it is an unconstrained least-squares problem, so it can be solved with standard optimization methods. Moreover, it only involves expectations over the reference and data distributions, and can thus be estimated directly from samples of $\alpha_0$ and $\alpha_1$ without requiring an explicit density estimation, which would be intractable in high dimensions.  
In practice, the vector field $v_t$ is parameterized by expressive neural architectures—most prominently U-Nets~\cite{ronneberger2015unet} and transformers~\cite{vaswani2017attention}—which have been successfully applied to images~\cite{ho2020denoising}, videos~\cite{ho2022video}, and other modalities~\cite{ramesh2022hierarchical}. One then minimizes~\eqref{eq:fm-linear-obj} with stochastic gradient descent. Sampling is performed by integrating $\dot x_t = v_t(x_t)$ from $t=0$ to $1$, starting at $x_0 \sim \alpha_0$.
This approach raises significant mathematical questions, particularly concerning the efficiency of these models. Key open issues, which we do not address in this section, include analyzing the sample complexity of score estimation~\cite{zhu2023sample}, quantifying the errors introduced by discretizing the dynamics~\cite{benton2023nearly}, and understanding how these errors propagate through the overall pipeline~\cite{hurault2025score}. A deeper theoretical understanding of these aspects is crucial to mitigating the high computational cost of diffusion models, which, despite their outstanding empirical performance, remain expensive due to the large number of integration steps they require.

\subsection{Special cases: Gaussian latent distributions} We now focus on the important case where $\alpha_0 = \mathcal{N}(0,\Sigma_0)$ is a Gaussian. 

\paragraph{Gaussian data.} 

To illustrate the method in a simple setting, the following proposition provides a closed-form expression for the flow-matching field when both the reference and data distributions are Gaussian.  

\begin{proposition}[\cite{hurault2025score}]\label{prop:gaussian-diffusion}
	If $\alpha_0 = \mathcal{N}(0,\Sigma_0)$ and $\alpha_1 = \mathcal{N}(0,\Sigma_1)$, then for $t \in [0,1]$ the interpolant is Gaussian $\alpha_t = \mathcal{N}(0,\Sigma_t)$ with covariance $\Sigma_t := (1-t)^2\Sigma_0 + t^2\Sigma_1$,
	and the associated velocity field is linear $v_t(z) = A_t z$ with
	$A_t = \bigl(-(1-t)\Sigma_0 + t \Sigma_1 \bigr)\Sigma_t^{-1}$.
\end{proposition}

Thus, in the Gaussian setting the velocity field $v_t$ is linear, and the transport map $T$ obtained by integrating $v_t$ is also linear.  
If the covariance matrices $\Sigma_0$ and $\Sigma_1$ do not commute, $A_t$ is not symmetric, and $v_t$ is not a gradient field. In this case, flow matching is not equivalent to the least-squares formulation~\eqref{eq:least-square-field}. However, when $\Sigma_0 = \mathrm{Id}_d$ (a common situation in practice), the flow-matching vector field $v_t$ is a gradient, and the map obtained by integrating it coincides with the optimal transport map described in Section~\ref{sec:ot}. Beyond this Gaussian case, the geometry and properties of flow matching remain largely unexplored; see~\cite{lavenant2022flow} for recent progress in this direction.

\paragraph{Connection with diffusion models}

If $\alpha_0$ is isotropic $\Sigma_0=\mathrm{Id}_d$, but $\alpha_1$ is arbitrary, the flow-matching method is equivalent, after an exponential reparameterization of time, to the diffusion model framework~\cite{ho2020denoising,song2020score}. A key mathematical advantage of this case is that the vector field $v_t$ can be expressed explicitly in terms of the interpolated law $\alpha_t$, or its density when it exists, as shown below.  

\begin{proposition}[Optimal vector field]\label{prop:flow}
When $\alpha_0 = \mathcal{N}(0,\mathrm{Id}_d)$, the vector field $v_t$ in~\eqref{eq:fm-linear-ce} is given by
\[
v_t(x) 
= \frac{1}{t}\,x \;+\; \frac{1-t}{t}\,\nabla\log \alpha_t(x),
\qquad (x\in\mathbb R^{d},\;t\in(0,1)).
\]
\end{proposition}

This result follows directly from Tweedie's identity~\cite{efron2011tweedie}, which relates the score function $\nabla\log \alpha_t$ to optimal Gaussian denoisers. It shows that in the diffusion model setting, $v_t$ is a gradient field, and therefore coincides with the solution of the constrained least-squares problem~\eqref{eq:least-square-field}. Extending this characterization to more general reference distributions $\alpha_0$ remains an open research question.

\section{Optimal Transport}
\label{sec:ot}

Instead of prescribing an explicit interpolant $(\alpha_t)_t$ as in the previous section, we now take a different perspective by seeking an “optimal” interpolation, in the sense first introduced by Monge in the 18th century~\cite{monge1781memoire}. This approach imposes a gradient structure on the displacement, opening the door to both efficient computational methods and applications beyond sampling.  
Monge's original formulation was a static problem of finding the most efficient rearrangement of mass. It was later reformulated by Kantorovich in the 20th century~\cite{kantorovich1942translation} in terms of couplings, which established the foundation of linear programming. The modern mathematical revival of optimal transport came with Brenier's theorem~\cite{brenier1991polar}, which showed that under mild assumptions the optimal transport map is the gradient of a convex function, thereby revealing its deep geometric structure.
Beyond mathematics, OT has found deep connections with statistics, where it provides sharp quantitative bounds for the convergence of empirical measures~\cite{dudley1969speed}. In machine learning, it has become a powerful tool to evaluate distances between distributions and to design generative models that learn from samples. A breakthrough came with the introduction of entropic regularization~\cite{cuturi2013sinkhorn}, which enabled scalable algorithms and inspired applications in generative adversarial networks (GANs) and other generative frameworks~\cite{genevay2018learning}.
A dynamic counterpart was introduced by Benamou and Brenier~\cite{benamou2000computational}, where OT is formulated as the minimization of kinetic energy along flows, leading to a convex optimization problem with wide-ranging applications.

\subsection{Static Optimal Transport}\label{subsec:static-ot}

Given two probability measures $\alpha_0,\alpha_1 \in \mathcal{P}_2(\mathbb{R}^d)$ (finite second moments), the \emph{Monge problem} seeks a measurable map $T:\mathbb{R}^d\to\mathbb{R}^d$ that transports $\alpha_0$ to $\alpha_1$,
$T_{\sharp}\alpha_0=\alpha_1$,  while minimizing the quadratic transportation cost:
\begin{equation}\label{eq:monge}
	\inf_{T_{\sharp}\alpha_0=\alpha_1} \int_{\mathbb{R}^d} \|x-T(x)\|^2 \,\mathrm{d}\alpha_0(x).
\end{equation}
This problem can be extended in many ways; in particular, the quadratic distance $\|x-T(x)\|^2$ can be replaced by more general costs~\cite{santambrogio2015optimal}.

Although this is not the focus of this review, problem~\eqref{eq:monge} is notoriously difficult to analyze theoretically and to solve numerically. For instance, an optimal map $T$ may fail to exist, and even the admissible set $\{T:T_{\sharp}\alpha_0=\alpha_1\}$ can be empty—for example, when $\alpha_0$ is discrete while $\alpha_1$ has a density.  
The canonical workaround is the convex linear programming relaxation introduced by Kantorovich~\cite{kantorovich1942translation}, which replaces maps $x \mapsto T(x)$ with couplings $\pi \in \Pi(\alpha_0,\alpha_1)$, i.e., probability measures on $\mathbb{R}^d \times \mathbb{R}^d$ with marginals $(\alpha_0,\alpha_1)$. The Kantorovich relaxation reads  
\begin{equation}\label{eq:kantorovich}
  \mathrm{W}_2(\alpha_0,\alpha_1)^2 \;=\; \inf_{\pi\in\Pi(\alpha_0,\alpha_1)} \iint \|x-y\|^2 \,\mathrm{d}\pi(x,y).
\end{equation}
When an optimal Monge map $T$ exists, it induces an optimal coupling $\pi^\star=(\mathrm{Id},T)_{\sharp}\alpha_0$, and the values of~\eqref{eq:monge} and~\eqref{eq:kantorovich} coincide.  
The square root of the optimal value, $\mathrm{W}_2(\alpha_0,\alpha_1)$, is the 2-Wasserstein distance. It defines a geometric notion of distance on the space of probability measures. As we will show next, this distance admits a rich geodesic structure and can be used to design optimization schemes directly over the space of probability distributions.
A cornerstone result in optimal transport is Brenier's theorem, which characterizes when the Monge problem~\eqref{eq:monge} admits a unique solution $T$ and describes its structure. In particular, it shows that $T$ is the gradient of a convex function, and hence a monotone map, in the sense that 
$\langle T(x)-T(y),\,x-y\rangle \;\ge\; 0$ whenever the representative is defined at $x$ and $y$.
In one dimension, $T$ reduces to the unique increasing map pushing $\alpha_0$ to $\alpha_1$. In higher dimensions, however, the situation is more subtle, since not every monotone map is the gradient of a convex function.  

\begin{theorem}[Brenier~\cite{brenier1991polar}]\label{thm:brenier}
Assume $\alpha_0$ is absolutely continuous with respect to Lebesgue measure on $\mathbb{R}^d$. Then there is a unique $T=\nabla \phi$ such that $T_\sharp \alpha_0=\alpha_1$ and this $T$ is the unique solution to the Monge problem~\eqref{eq:monge}.
\end{theorem}

\subsection{Dynamical Optimal Transport}

The Monge problem~\eqref{eq:monge} seeks an instantaneous displacement from $t=0$ to $t=1$ through a transport map $T$. An alternative viewpoint is to compute an interpolation $t \mapsto \alpha_t$ that minimizes the least-squares energy~\eqref{eq:least-square-field}. Unlike the flow-matching setting studied in Section~\ref{sec:genmodels}, where the interpolation $(\alpha_t)_t$ is prescribed, here we only assume knowledge of the initial and final distributions $\alpha_0$ and $\alpha_1$.  
A fundamental result, proved by Benamou and Brenier, is that the minimal value of this ``geodesic'' energy coincides with the squared Wasserstein-2 distance. Moreover, the optimal pair $(v_t,\alpha_t)$ solving this variational problem can be constructed by evolving particles along straight lines driven by the Monge map $T$ that solves the original Monge problem.

\begin{theorem}[Benamou--Brenier~\cite{benamou2000computational}]
For probability measures $\alpha_0,\alpha_1$ with finite second moments,
\begin{equation}\label{eq:benamou-brenier}
  \mathrm{W}_2^2(\alpha_0,\alpha_1)
    = \inf_{(\alpha_t,v_t)} \int_0^1 \int_{\mathbb{R}^d} |v_t(x)|^2\, \d\alpha_t(x)\, \d t,
\end{equation}
where the infimum is over $(\alpha_t,v_t)$ solving
$\partial_t\alpha_t+\mathrm{div}(\alpha_t v_t)=0$ with $\alpha_{t=0}=\alpha_0$ and $\alpha_{t=1}=\alpha_1$.  
If $\alpha_0$ has a density and $T$ is the optimal Monge map $T_{\sharp}\alpha_0=\alpha_1$, the minimizer is
\begin{equation}\label{eq:static-to-dynamic}
  \alpha_t=((1-t)\operatorname{Id}+tT)_{\sharp}\alpha_0,
  \qquad
  v_t\bigl((1-t)x+tT(x)\bigr)=T(x)-x,\quad t\in[0,1].
\end{equation}
\end{theorem}

Theorem~\ref{thm:brenier} shows that the optimal displacement $T$ has a gradient form. The dynamic problem~\eqref{eq:benamou-brenier}, assuming the optimal path $(\alpha_t)_t$ is known, reduces to a least-squares problem of the form~\eqref{eq:least-square-field}. Following the discussion of Section~\ref{sec:least-square-field}, this implies that the associated velocity field $v_t$ is itself a gradient field.  
The interpolation formula~\eqref{eq:static-to-dynamic} bridges the static formulation~\eqref{eq:monge} and the dynamic one~\eqref{eq:benamou-brenier}. It shows that once the optimal map $T$ is known, the interpolating curve $(\alpha_t)_t$ is obtained by moving particles along straight lines from $X_0$ to $T(X_0)$. This should not be confused with the stochastic interpolation described in Section~\ref{sec:stoch-int}, where $\alpha_t$ is defined by connecting \emph{independent} random vectors $X_0 \sim \alpha_0$ and $X_1 \sim \alpha_1$. In optimal transport, by contrast, one connects $X_0$ with its coupled image $T(X_0)$, which are not independent. This subtle but fundamental difference explains why OT interpolations and flow matching behave so differently: flow matching proceeds through convolution, producing smoother interpolations.  
Finally, while not the focus of this survey, it is worth noting that the dynamic problem~\eqref{eq:benamou-brenier} is non-convex in $(\alpha_t,v_t)_t$. However, after the change of variables $m_t := \alpha_t v_t$, it becomes convex in $(m_t,\alpha_t)_t$. This elegant property allows geodesic interpolation to be solved efficiently using convex optimization techniques.

\paragraph{Gaussian case.}  

In the Gaussian setting, the Monge problem admits a closed-form solution: the optimal transport map is linear, and the Wasserstein-2 distance can be computed explicitly.  

\begin{proposition}\label{prop:wass-gaussian}
For $\alpha_0=\mathcal{N}(0,\Sigma_0)$ and $\alpha_1=\mathcal{N}(0,\Sigma_1)$, the optimal transport map is linear 
    \[
    T(x) \;=\; \Sigma_0^{-\tfrac{1}{2}}
      \Bigl(\Sigma_0^{\tfrac{1}{2}}\Sigma_1\Sigma_0^{\tfrac{1}{2}}\Bigr)^{\tfrac{1}{2}}
      \Sigma_0^{-\tfrac{1}{2}} \, x,
    \]
    and the squared Wasserstein-2 distance is the so-called squared Bures distance~\cite{dowson1982frechet,takatsu2011wasserstein}  
    \[
    \mathrm{W}_2^2(\alpha_0,\alpha_1) \;=\; 
    \mathrm{Tr}\!\Bigl(\Sigma_0+\Sigma_1 - 2\bigl(\Sigma_0^{1/2}\Sigma_1\Sigma_0^{1/2}\bigr)^{1/2}\Bigr).
    \]
    The geodesic interpolation $\alpha_t=((1-t)\mathrm{Id}+tT)_\sharp \alpha_0$ is Gaussian with covariance $\Sigma_t = \bigl((1-t)\mathrm{Id}+tT\bigr)\,\Sigma_0\,\bigl((1-t)\mathrm{Id}+tT\bigr)^{\top}$.
\end{proposition}

This proposition shows that, in general, the optimal transport map and the map obtained by integrating the flow-matching field in Proposition~\ref{prop:gaussian-diffusion} are different. They coincide only when $\Sigma_0$ and $\Sigma_1$ commute, which in particular holds when $\Sigma_0=\mathrm{Id}_d$ (see Proposition~\ref{prop:flow}). This highlights the structural gap between flow-matching and optimal transport interpolations.

\section{Wasserstein Gradient Flows}
\label{sec:wass-flow}

Perhaps the most influential development of this idea of dynamical OT is the notion of Wasserstein gradient flows, where flows are steered by the gradient of an energy functional. 
Beyond interpolating between two probability distributions, one can then study the evolution of measures that minimize a functional $f : \Pp(\RR^d) \to \RR \cup \{+\infty\}$. 
This idea was pioneered in the seminal JKO scheme~\cite{jordan1998variational} and further developed in the comprehensive theory of gradient flows in metric spaces by Ambrosio, Gigli, and Savaré~\cite{ambrosio2008gradient}. This framework allows nonlinear PDEs such as the porous medium equation~\cite{otto2001geometry} or crowd-motion models with congestion~\cite{maury2010macroscopic} to be interpreted as gradient flows in Wasserstein space, yielding new insights into their well-posedness and long-time behavior.
In machine learning, Wasserstein gradient flow techniques have been used to analyze optimization dynamics. A notable line of work views the training of shallow neural networks as a gradient flow over probability measures on the space of neurons. This approach was initiated by Chizat and Bach~\cite{chizat2018global} and extended in subsequent works~\cite{mei2018meanfield,rotskoff2018parameters,wei2019regularization}, showing convergence properties and connections with mean-field limits.

\subsection{The JKO scheme}
\label{sec:jko}

The starting point is the minimizing movement scheme, also known as the JKO scheme~\cite{jordan1998variational}. Given a time step $\tau>0$, one defines the discrete evolution
\begin{equation}\label{eq:jko-discr}
    \alpha_{t+\tau} := \arg\min_{\alpha \in \Pp(\RR^d)}
    \left\{ \frac{1}{2\tau} \mathrm{W}_2^2(\alpha_t,\alpha) + f(\alpha) \right\}.
\end{equation}
This is the analogue in Wasserstein space of the implicit Euler method for gradient flows in Euclidean space.  
To see this, suppose $\alpha_t=\delta_{x(t)}$ is a Dirac mass and set $h(x)=f(\delta_x)$, assuming $f$ can be evaluated on Dirac masses. Then~\eqref{eq:jko-discr} reduces to
\begin{equation}\label{eq:prox-particles}
    x(t+\tau) = \arg\min_x \left\{ \frac{1}{2\tau}\|x-x(t)\|^2 + h(x) \right\},
\end{equation}
which is precisely the implicit Euler step for the ODE $\dot x(t) = -\nabla h(x(t))$. As $\tau \to 0$, the JKO iteration~\eqref{eq:jko-discr} converges (under suitable conditions) to a continuous-time evolution $(\alpha_t)_{t\geq 0}$, which can be described by the continuity PDE~\eqref{eq:eulerian-advection}.
The advantage of this approach is that one can single out a vector field $v_t$ with a closed-form expression, the negative Wasserstein gradient $v_t = -\nabla_{\mathrm{W}} f(\alpha_t)$, which is formally a gradient field (in line with the study of Section~\ref{sec:least-square-field}), defined as
\[
    \nabla_{\mathrm{W}} f(\alpha) = \nabla [\delta f(\alpha)], 
\]
see~\cite{ambrosio2008gradient}. Here $\delta f(\alpha)$ is the first variation of $f$, i.e., a function in $\mathcal{C}(\RR^d)$ satisfying
\[
    f((1-\tau)\alpha+\tau\beta) 
    = f(\alpha) + \tau \int \delta f(\alpha)(x)\,\d(\beta-\alpha)(x) + o(\tau).
\]
Thus the Wasserstein gradient flow is driven by the potential $\delta f(\alpha)$.
The Wasserstein gradient flow limiting PDE thus reads
\begin{equation}\label{eq:wassflow-pde}
    \partial_t \alpha_t + \mathrm{div}(-\nabla_{\mathrm{W}} f(\alpha_t)\, \alpha_t ) = 0.
\end{equation}
We now detail several special cases of interest, both for the study of partial differential equations and for the training of neural networks.

\subsection{Diffusions}

A striking illustration of this framework is the following result.

\begin{theorem}[Jordan–Kinderlehrer–Otto~\cite{jordan1998variational}]\label{thm:jko-heat}
Let $f(\alpha) = \int_{\RR^d} \log\!\left(\frac{\d \alpha}{\d x}(x)\right)\, \d \alpha(x)$
be the Shannon–Boltzmann negative entropy (defined for absolutely continuous $\alpha$). Then the Wasserstein gradient flow of $f$ is the heat equation
$\partial_t \alpha_t = \Delta \alpha_t$. 
\end{theorem}

This identification provides a new variational interpretation of diffusion. Traditionally, the heat equation arises as the $L^2$ gradient flow of the Dirichlet energy
$f(\alpha) = \int_{\RR^d} \|\nabla \rho(x)\|^2 \,\d x$ where $\rho = \frac{\d \alpha}{\d x}$, with respect to the $L^2$ Euclidean distance between densities. The JKO formulation instead replaces the geometry (Euclidean $\to$ Wasserstein) and the functional (Dirichlet $\to$ negative entropy), yet recovers the same PDE.  
From the transport viewpoint, this means diffusion can be understood as the motion of particles seeking to spread mass as uniformly as possible—that is, as the gradient flow that maximizes entropy.

Other entropy functionals lead to nonlinear diffusion equations. For example, generalized entropies of the form 
\begin{equation}\label{eq:gen-entropies}
       f(\alpha) = \int g\!\left(\frac{\mathrm{d} \alpha}{\mathrm{d} x}\right) \mathrm{d} x,
\end{equation}
for a one-dimensional function $g(s)$, give rise to nonlinear diffusions of the type
$\frac{\partial \alpha_t}{\partial t} = \Delta(\tilde{g}(\alpha_t))$,
where $\tilde{g}$ is defined, up to an additive constant, by $\tilde{g}'(s)=s g''(s)$. For instance, $g(s) = s \log(s)$ corresponds to Theorem~\ref{thm:jko-heat}, while $g(s) = s^p/(p-1)$ with $p > 1$ yields the porous medium (slow diffusion) equation~\cite{otto2001geometry}.
This interpretation of diffusion as a Wasserstein gradient flow is crucial for analyzing long-time behavior and convergence to stationary equilibria.
A celebrated and non-trivial theorem by McCann~\cite{mccann1997convexity} gives geodesic convexity for functionals of the form~\eqref{eq:gen-entropies}, with $g : \RR^+ \to \RR$, under the usual displacement-convexity conditions: $g(0)=0$, $g$ is convex and superlinear, and the map $s \mapsto g(s^{-d}) s^d$ is convex and decreasing.
In particular, if $(\alpha_t)_{t\in[0,1]}$ is a Benamou--Brenier geodesic solving~\eqref{eq:benamou-brenier}, then the function $t \mapsto f(\alpha_t)$ is convex. This convexity is a key ingredient in convergence analyses for the PDE~\eqref{eq:wassflow-pde}; quantitative or linear rates require stronger assumptions, such as suitable confinement or functional inequalities.
Examples of such $g$ include $g(s)=s^q$ for $q>1$ and the Shannon entropy $g(s)=s \log(s)$.
These convexity arguments can be extended to handle mild non-convexity and have been applied in machine learning to study the convergence to equilibrium of noisy gradient descent~\cite{mei2018meanfield,zhang2022meanfield}, which is often used as a heuristic to understand the behavior of stochastic gradient descent in training.

\subsection{Interacting particles}
\label{sec-interparticles}

Diffusion flows smooth out densities and therefore do not preserve discrete distributions. While this property is useful for mathematical analysis of convergence, for instance through geodesic convexity arguments, it is less suitable for modeling certain machine learning problems. A prominent example is the training of neural networks without artificial noise injection, which is better described as the evolution of large but finite particle systems, as discussed in Section~\ref{sec:mlp}.  
If $f(\alpha)$ can be evaluated on discrete measures and $\nabla_{\mathrm{W}} f(\alpha)$ defines a regular vector field, then the flow \eqref{eq:wassflow-pde} preserves the discrete structure
$\alpha_t = \frac{1}{n} \sum_{i=1}^n \delta_{x_i(t)}$, so that the particles $X(t) := (x_i(t))_{i=1}^n$ evolve according to a system of coupled ODEs
\begin{equation}\label{eq:wassflows-particles}
    \dot{x}_i(t) = -[\nabla_{\mathrm{W}} f(\alpha_t)](x_i(t)), \qquad
    \alpha_t=\frac1n\sum_{i=1}^n \delta_{x_i(t)}.
\end{equation}
Equivalently, if $F(X) := f( \tfrac{1}{n} \sum_i \delta_{x_i})$ is viewed as a function of the particle positions with the standard Euclidean metric, then $\dot X(t)=-n\nabla F(X(t))$; the factor $n$ is only a convention coming from the normalization of the empirical measure.
A prototypical example is provided by quadratic interaction functionals of the form
\begin{equation}\label{eq:quadratic-func}
    f(\alpha) := \iint k(x,y)\,\d\alpha(x)\,\d\alpha(y).
\end{equation}
When $k(x,y) = h(x)$ does not depend on $y$, the Wasserstein gradient reduces to $\nabla_{\mathrm{W}} f(\alpha) = \nabla h$, independently of $\alpha$, and the particles move independently according to $\dot{x} = -\nabla h(x)$, which is exactly the setup considered in~\eqref{eq:prox-particles}. In the general case, however, the Wasserstein gradient induces couplings between particles (here we assume symmetry of $k$ for simplicity),
 $\nabla_{\mathrm{W}} f(\alpha)(x) = 2 \int \nabla_x k(x,y)\,\d\alpha(y)$, 
which leads to the particle system
\[
    \dot{x}_i(t) = -\frac{2}{n} \sum_{j=1}^n \nabla_x k(x_i(t),x_j(t)).
\]
Interacting particle models of this type are ubiquitous in physics, where they describe aggregation and self-organization phenomena~\cite{yao2014asymptotic}, in biology, where they are used to capture collective dynamics such as swarming~\cite{carrillo2019particle}, and in applied mathematics, where they arise in the study of mean-field limits and kinetic equations~\cite{carrillo2020derivation}. In general, however, quadratic interaction energies of the form~\eqref{eq:quadratic-func} are not geodesically convex, which makes the analysis of convergence to equilibrium particularly challenging.

\subsection{Trajectory inference for genomics}

Beyond the theoretical analysis of MLP training, another major application of Wasserstein flows has emerged in computational biology. In single-cell genomics, OT-based methods are now central to modeling how cell populations evolve over time, providing a mathematical framework to reconstruct developmental trajectories, differentiation pathways, and multi-omics dynamics from high-dimensional data.
Single-cell genomics~\cite{tanay2017scaling}—including scRNA-seq, scATAC-seq, CITE-seq, and related assays—has transformed biology by making it possible to measure heterogeneous cell populations at the resolution of individual cells, revealing continuous cell-state landscapes and rare subpopulations that bulk methods blur. Formally, such data are naturally represented as a probability measure $\alpha$ on an omics space $\mathcal{X}$. For scRNA-seq, $\mathcal{X}$ can be taken as a (typically high-dimensional) gene-expression space, where each cell corresponds to a noisy point sample and $\alpha$ describes the empirical distribution over that space. This level of detail enables the study of dynamic phenomena—e.g., differentiation during early development or the response to cancer therapy—through time-resolved snapshots. A central problem is therefore: given observations of snapshot measures $\alpha_t$ at a few time points $t$ (for example, across biopsies), recover the velocity field $v_t$ governing population evolution. This is precisely the trajectory-inference setting, and many of the methods discussed here are being explored for this purpose.
Inferring cell-fate trajectories from spatially resolved profiles collected at multiple times has thus become a key objective, motivating new computational strategies. Wasserstein gradient-flow learning offers a promising direction for longitudinal sequencing data, where a neural network encodes a potential landscape of differentiation over time.
Trajectory inference methods tailored to single-cell data~\cite{saelens2019comparison} aim to reconstruct developmental paths and branching structures. Waddington-OT models trajectories by estimating probabilistic cell–cell transitions between successive time points~\cite{schiebinger2019optimal}. Another family of OT-based approaches posits a continuous population dynamics model by training neural networks to represent a generalized notion of velocity~\cite{tong2020trajectorynet}. A further OT-based line of work seeks to learn a potential function $f(\alpha)$ that drives a causal differentiation model, closely aligned with the Wasserstein gradient-flow formalism described in Section~\ref{sec:jko}~\cite{hashimoto2016learning,yeo2021generative,bunne2022proximal}. Casting differentiation as the minimization of a potential has deep roots in systems biology and provides a mathematical expression of Waddington's epigenetic landscape~\cite{allen2015compelled}.
Finally, integrating measurements from multiple modalities (multi-omics)—for example, combining RNA-seq with spatial cell coordinates—poses substantial challenges that can be addressed using extensions of optimal transport~\cite{shen2025inferring,huizing2024stories}.

\section{Training Neural Networks}
\label{sec:mlp}

In this section, we show how the Wasserstein flow framework from the previous part applies to the analysis of shallow network training dynamics, where the network is represented by the probability distribution of its weights, and this distribution evolves by minimizing the training error.

\subsection{Training as Wasserstein Flows}

A basic but expressive example of neural architectures is the two-layer multilayer perceptron (MLP), introduced in early work on feedforward networks~\cite{cybenko1989approximation,hornik1991approximation}. Such networks define functions $G_X : \RR^d \to \RR^{d'}$ of the form
\begin{equation}\label{eq:2layermlp}
    G_X(z) := \frac{1}{n} \sum_{i=1}^n v_i \, \sigma(\langle z, u_i \rangle),
\end{equation}
where each neuron is parameterized by an inner weight $u_i \in \RR^d$ and an outer weight $v_i \in \RR^{d'}$, collected into $x_i := (u_i,v_i)$ and $X=(x_i)_i$. Here $z \in \RR^d$ is the input, and $\sigma:\RR \to \RR$ is the activation function. The nonlinearity of $\sigma$ (it must be non-polynomial) ensures expressivity; the most widely used choice in modern practice is the rectified linear unit (ReLU), $\sigma(s) = \max\{s,0\}$~\cite{nair2010rectified}. Universal approximation results show that as $n\to\infty$, such models can approximate uniformly any continuous function on compact sets~\cite{cybenko1989approximation,hornik1991approximation}.  

Training amounts to adjusting the parameters $X=(x_i)_{i=1}^n$ using a dataset $(z_k,y_k)_{k=1}^N \subset \RR^d \times \RR^{d'}$ of input–output pairs. The standard approach is empirical risk minimization. For simplicity, we consider the quadratic loss,
\begin{equation}\label{eq:empirical-risk}
    F(X) := \frac{1}{N} \sum_{k=1}^N \| G_X(z_k) - y_k \|^2,
\end{equation}
and perform gradient descent (in practice, stochastic gradient descent and its variants~\cite{robbins1951sgd,bottou2018optimization}):
\[
    X_{t+\tau} = X_t - \tau \nabla F(X_t).
\]
As $\tau \to 0$, this iteration converges formally to continuous-time gradient descent; up to the normalization factor discussed in~\eqref{eq:wassflows-particles}, this is the particle form of the Wasserstein flow. Global convergence for large networks remains poorly understood, but significant progress has been made in the case of two-layer MLPs~\eqref{eq:2layermlp} by analyzing their mean-field limit.  

The key observation is that $G_X$ is invariant under permutation of neurons, so that it is natural to represent the weights by a probability measure $\alpha$ on $\RR^{d+d'}$. This leads to the mean-field formulation
\begin{equation}\label{eq:mlp-meanfield}
    g_\alpha(z) := \int_{\RR^{d+d'}} v\,\sigma(\langle z, u\rangle)\, \d\alpha(u,v),
    \qquad
    f(\alpha) := \frac{1}{N} \sum_{k=1}^N \| g_\alpha(z_k) - y_k \|^2.
\end{equation}
When $\alpha = \tfrac{1}{n}\sum_{i=1}^n \delta_{(u_i,v_i)}$, one recovers $g_\alpha = G_X$. The functional $f(\alpha)$ is therefore a continuous relaxation of~\eqref{eq:empirical-risk} valid even when $\alpha$ has a density. Minimizing $f(\alpha)$ via Wasserstein gradient flow corresponds to the particle dynamics of Section~\ref{sec-interparticles}. For quadratic loss, $f$ can be written, up to an additive constant, as a quadratic functional of the form~\eqref{eq:quadratic-func}, with an interaction kernel
\begin{equation}\label{eq:kernel-mlp}
    k\bigl((u,v),(u',v')\bigr)
    = \frac{1}{N}\sum_{k=1}^N \bigl\langle v,v' \bigr\rangle
        \sigma(\langle u,z_k\rangle)\sigma(\langle u',z_k\rangle)
      - \frac{1}{N}\sum_{k=1}^N
      \Bigl(
      \langle v, y_k \rangle \sigma(\langle u,z_k\rangle)
      + \langle v', y_k \rangle \sigma(\langle u',z_k\rangle)
      \Bigr).
\end{equation}
In general, global convergence of Wasserstein flows for quadratic interaction energies is out of reach. However, a breakthrough result by Chizat and Bach~\cite{chizat2018global} establishes convergence for this particular setting.  

\begin{theorem}[Chizat--Bach]\label{thm:chizat-bach}
Consider the two-layer MLP~\eqref{eq:mlp-meanfield} with ReLU activation. If the initialization $\alpha_0$ has a density with respect to Lebesgue measure on $\RR^{d+d'}$, then every convergent limit point of the Wasserstein gradient flow of $f(\alpha)$ is a global minimizer of $f$.
\end{theorem}

A related result holds for smooth bounded activations~\cite{chizat2018global}, and was later refined to properly handle the non-smoothness of ReLU~\cite{wojtowytsch2020convergence}. While these theorems assume an infinite number of neurons at initialization, Chizat and Bach also show that the discrete particle system with finitely many neurons approximates the mean-field limit over controlled time horizons. Since two-layer MLPs are universal approximators, the infimum of the empirical risk can be made arbitrarily small, and it can be zero under standard realizability or interpolation assumptions.
The key remaining challenge is that these results are non-quantitative: it remains unclear how large $n$ must be in order to guarantee global convergence, even under suitable assumptions on the data distribution $(z_i,y_i)_i$.

\subsection{Linear neural network under a Gaussian neuron distribution}

In the spirit of Propositions~\ref{prop:gaussian-diffusion} and~\ref{prop:wass-gaussian}, we consider a Gaussian distribution of neurons $\alpha=\mathcal N(0,\Sigma)$ over $(u,v)\in\RR^d\times\RR^{d'}$, with block covariance
\[
\Sigma \;=\;
\begin{pmatrix}
\Sigma_{UU} & \Sigma_{UV}\\[2pt]
\Sigma_{VU} & \Sigma_{VV}
\end{pmatrix},
\qquad
\Sigma_{UV} \;=\; \int u v^\top\,\d\alpha(u,v),\quad
\Sigma_{VU}=\Sigma_{UV}^\top.
\]
We study the simplified linear network $\sigma(s)=s$, a classical toy model for training dynamics~\cite{saxe2013exact,gunasekar2018implicit}. In this case,
\[
g_\alpha(z) \;=\; \int v\,\langle z,u\rangle\,\d\alpha(u,v)
\;=\; \Sigma_{VU}\,z,
\]
so the predictor is linear with weight matrix $\Sigma_{VU} \in \RR^{d'\times d}$. Let the data moments be
\[
S_{zz}:=\frac1N\sum_{k=1}^N z_k z_k^\top \in\RR^{d\times d},\qquad
S_{yz}:=\frac1N\sum_{k=1}^N y_k z_k^\top \in\RR^{d'\times d}.
\]
Then the empirical risk from~\eqref{eq:mlp-meanfield} becomes
\[
f(\alpha)\;=\;\frac1N\sum_{k=1}^N \|\Sigma_{VU} z_k - y_k\|^2
\;=\;\operatorname{Tr}\!\big(\Sigma_{VU} S_{zz} \Sigma_{UV}\big)\;-\;2\,\operatorname{Tr}\!\big(S_{yz}\,\Sigma_{UV}\big)\;+\;\text{const},
\]
which depends on $\alpha$ only through $\Sigma_{UV}$ and $\Sigma_{VU}$. The following proposition shows that the training dynamics, given by the Wasserstein flow~\eqref{eq:wassflow-pde}, preserve Gaussian neuron distributions.

\begin{proposition}[Wasserstein flow for linear networks closes on Gaussians]\label{prop:gaussian-mlp}
Assume $S_{zz}$ is full rank and $\sigma(s)=s$. The Wasserstein gradient-flow velocity field $-\nabla_{\mathrm{W}} f(\alpha)$ is the linear vector field
\[
[-\nabla_{\mathrm{W}} f(\alpha)](u,v) \;=\; -2 \Bigl( (S_{zz}\Sigma_{UV}-S_{zy})v,\; (\Sigma_{VU}S_{zz}-S_{yz})u\Bigr),
\]
with $S_{zy}=S_{yz}^\top$. Consequently, if $\alpha_{t=0}=\mathcal N(0,\Sigma_0)$, then $\alpha_t=\mathcal N(0,\Sigma_t)$ for all $t$, where the covariance satisfies the matrix ODE
\begin{equation}\label{eq:gaussian-mlp}
\dot\Sigma_t \;=\; B_t \Sigma_t + \Sigma_t B_t^\top,
\qquad
B_t \;:=\;
-2 
\begin{pmatrix}
0 & (S_{zz}\Sigma_{UV,t}-S_{zy})\\[2pt]
(\Sigma_{VU,t} S_{zz}-S_{yz}) & 0
\end{pmatrix}.
\end{equation}
\end{proposition}

This shows that the training dynamics reduce to the Riccati-type ODE~\eqref{eq:gaussian-mlp} over the space of covariance matrices~\cite{abou2003matrix}. The system is locally well posed, and the covariance description remains valid as long as the solution stays finite. When it converges to a finite minimizer, the limiting predictor is characterized by $\Sigma_{VU,\infty} = S_{yz}\,S_{zz}^{-1}$, or equivalently $\Sigma_{UV,\infty} = S_{zz}^{-1}S_{zy}$, so that $g_\infty$ coincides with the least-squares linear estimator.

\subsection{Toward deeper networks}

The convergence results discussed above apply only to shallow two-layer networks in a mean-field regime. Extending such guarantees to deeper networks remains a largely open problem. From an empirical standpoint, it is widely observed that very deep and large neural networks can be trained effectively despite their non-convex loss landscapes. However, a rigorous mathematical explanation of this phenomenon is still missing. Importantly, these are optimization results: while a network can often be trained to a global minimizer, not all global minima generalize equally well to unseen data.  
Most theoretical analyses (see, e.g.,~\cite{jacot2018neural,lee2019wide,arora2019fine}) rely on a suitable rescaling of network parameters in the infinite-width limit, leading to the so-called \emph{kernel regime}. In this regime, training dynamics are approximated by a quadratic optimization problem involving the \emph{neural tangent kernel} (NTK)~\cite{jacot2018neural}. While powerful for analysis, this approximation poorly reflects real-world training, since it predicts a ``lazy training'' behavior associated with limited feature learning and poor generalization~\cite{chizat2019lazy}.  

To capture richer dynamics, it is necessary to account for \emph{skip connections}, as introduced in residual networks (ResNets)~\cite{he2016deep}. With $S$ residual blocks, one can model a block as a two-layer MLP correction with scaling $1/S$:  
\begin{equation}\label{eq:skip-connection}
    z_{\ell+1} := z_\ell + \frac{1}{S}\, g_{\alpha_\ell}(z_\ell),
    \qquad \ell=0,\ldots,S-1,
\end{equation}
where $z_\ell$ denotes the activations at layer $\ell$ (with $z_0$ the input).  
Taking the infinite-depth limit $S\to\infty$ yields the \emph{neural ODE} model~\cite{chen2018neural,weinan2017proposal,haber2017stable},  
\begin{equation}\label{eq:neuralode}
    \frac{\d z}{\d s}(s) = g_{\alpha(s)}(z(s)), \quad s\in[0,1],
\end{equation}
parameterized by a family $\Alpha=(\alpha(s))_{s\in[0,1]}$. The network output is defined as $G_\Alpha(z(0)) := z(1)$. Training amounts to minimizing the empirical risk $f(\Alpha) = \tfrac{1}{N}\sum_{k=1}^N \|G_\Alpha(z_k)-y_k\|^2$ via gradient descent in $\Alpha$.  

As shown in~\cite{barboni2024understanding}, the gradient descent dynamics for the ResNet~\eqref{eq:skip-connection} converge, as $S\to\infty$, to a Wasserstein gradient flow with respect to a generalized metric over depth-dependent distributions:  
\[
    \mathcal{W}_2(\Alpha,\Alpha')^2 := \int_0^1 \mathrm{W}_2(\alpha(s),\alpha'(s))^2\,\d s.
\]
This distance can be seen as a \emph{conditional optimal transport} metric over paths of measures. Using this formalism, Barboni et al.~\cite{barboni2024understanding} prove local convergence toward global minimizers, despite the strong non-convexity of the loss function. A key point is that the loss satisfies a Polyak–Łojasiewicz (PL) inequality~\cite{karimi2016linear} near minimizers, ensuring that the non-convexity is not pathological.  
Extending these local results to global convergence remains an open problem. One major challenge is that the loss functional $f(\Alpha)$ is not coercive on the path space endowed with $\mathcal{W}_2$: under certain initializations, parameters may degenerate and diverge to infinity. Understanding how to control such degeneracies, and whether realistic initialization schemes avoid them, is an important direction for future work.

\section{Transformers}
\label{sec:transformers}

Transformers, first introduced in~\cite{vaswani2017attention}, have become one of the most influential architectures in deep learning, powering state-of-the-art results across natural language processing, vision, and beyond. A common paradigm is to \emph{tokenize} data~\cite{lewis2020pretraining}: split it into subparts called tokens, represented as vectors $(x_i)_{i=1}^n$ in a latent space $\RR^d$. For text, tokens correspond to subword units enriched with positional embeddings to encode ordering in the sentence~\cite{sennrich2016neural,vaswani2017attention}; for vision, tokens are typically patches of pixels mapped into embeddings~\cite{dosovitskiy2021an}. The network then processes this point cloud of tokens through a repeated composition of three operations: normalization layers, multilayer perceptrons (MLPs), and, most importantly, self-attention blocks. 
While normalization and MLP layers act independently on each token, self-attention couples tokens together, enabling the capture of arbitrarily long-range dependencies and structural relations. Note that MLP layers hold most of the network's parameters in large models~\cite{kaplan2020scaling}.

The following sections present a mathematical model for deep transformers, which follows from the observation that transformers are naturally permutation-equivariant when applied to sets of tokens without causal constraints (e.g., in vision tasks~\cite{dosovitskiy2021an}), which motivates viewing the token collection as a probability distribution $\alpha$ on the token space. This suggests a probabilistic and measure-theoretic lens for analysis. 
Deep transformers, which stack many layers~\cite{lin2022survey}, can be modeled in the infinite-depth limit as continuous flows $t \mapsto \alpha_t$, akin to neural ODEs~\cite{chen2018neural,weinan2017proposal,haber2017stable}. This leads to a conservation law or PDE governing the evolution of $\alpha_t$ through depth, formalized in~\cite{sander2022sinkformers}. Recent work has explored the long-time behavior of such ``transformer PDEs'', including clustering phenomena on the sphere~\cite{zhong2022neural,geshkovski2023emergence,geshkovski2023mathematical} and tractable cases such as Gaussian mixtures~\cite{castin2025unified}. However, because the attention map is neither smooth nor linear in the underlying measure~\cite{castin2023smooth}, the analysis is delicate: the transformer PDE is more closely related to Vlasov- or aggregation-type equations~\cite{vlasov1968vibrational,Go03} than to classical gradient flows. In fact, it does not generally admit a Wasserstein gradient flow structure~\cite{sander2022sinkformers}, although recent work suggests that suitable modifications of the metric can restore one~\cite{burger2025analysismeanfieldmodelsarising}.

\subsection{Infinite number of tokens and mean-field limit}

An attention layer with skip connection and rescaling by $1/T$ (where $T$ is the network depth, i.e., the total number of layers in the transformer) defines the following update rule for the tokens:
\begin{equation}\label{eq:attention-discr}
    x_i \in \RR^d \;\mapsto\; x_i + \frac{1}{T} \sum_j \frac{\exp\langle Q x_i, K x_j \rangle \, V x_j}{\sum_{\ell} \exp\langle Q x_i, K x_\ell \rangle},
\end{equation}
where $\theta = (K, Q, V)$ collects the key, query, and value matrices, each in $\RR^{d \times d}$. This is a simplified version of the full multi-head attention mechanism, which combines several such maps and typically uses rectangular $(Q,K,V)$ matrices to reduce the number of parameters.

Transformers are designed to operate on a variable, often large, number $n$ of tokens during training and inference. It is therefore natural to analyze their behavior in the mean-field regime $n \to \infty$.  
Following Sander et al.~\cite{sander2022sinkformers}, we introduce the empirical measure of tokens
$\alpha = \frac{1}{n} \sum_{i=1}^n \delta_{x_i}$, 
and rewrite the transformer update in terms of $\alpha$ as
\begin{equation}\label{eq:transform-step}
    x_i \;\mapsto\; x_i + \frac{1}{T} \, \Gamma_\theta[\alpha](x_i),
    \qquad
    \Gamma_\theta[\alpha](x) := 
    \frac{\int e^{\langle Qx, Ky\rangle} V y \, \mathrm{d}\alpha(y)}{\int e^{\langle Qx, Kz\rangle} \, \mathrm{d}\alpha(z)}.
\end{equation}
For a discrete measure $\alpha = \tfrac{1}{n} \sum_i \delta_{x_i}$, this formulation reduces to the finite-token expression~\eqref{eq:attention-discr}.

\subsection{Infinitely deep transformer and the transformer PDE}

As $T\to\infty$, the discrete update in~\eqref{eq:transform-step} becomes an explicit Euler discretization of a coupled ODE system for the particles $(x_i(t))_{i=1}^n$ indexed by a continuous depth variable $t\in[0,1]$:
\begin{equation}\label{eq:transform-ode-coupled}
    \dot x_i(t) \;=\; \Gamma_{\theta_t}[\alpha_t]\bigl(x_i(t)\bigr),
    \qquad
    \alpha_t := \frac{1}{n}\sum_{i=1}^n \delta_{x_i(t)}.
\end{equation}
This is the neural-ODE viewpoint~\cite{chen2018neural,weinan2017proposal,haber2017stable} specialized to transformers. Equivalently, the time-evolving empirical measure $(\alpha_t)_{t\in[0,1]}$ solves the continuity equation (``transformer PDE'') introduced in~\cite{sander2022sinkformers}:
\begin{equation}\label{eq:transformer-pde}
    \partial_t \alpha_t \;+\; \nabla\!\cdot\!\bigl(\alpha_t \,\Gamma_{\theta_t}[\alpha_t]\bigr) \;=\; 0.
\end{equation}
Note that for simplicity, we omit the MLP and normalization layers in the transformer and assume a single head per attention layer.
Assume $\theta_t\equiv\theta=(Q,K,V)$ is fixed in time and suppose $V=Q^{\top}K$. Define the log-partition functional
\[
    \phi[\alpha](x) \;:=\; \log \int e^{\langle Qx,Ky\rangle}\,\mathrm{d}\alpha(y).
\]
A direct calculation shows that if one assumes $V=Q^{\top}K$, then $\Gamma_{\theta}[\alpha] =  \nabla_x \phi[\alpha]$. Hence in this case, the drift is a gradient field $v_\alpha=\nabla \phi[\alpha]$. 
However, as observed in~\cite{sander2022sinkformers}, $\phi[\alpha]$ is \emph{not} a first variation $\delta f(\alpha)$ of a functional on $\Pp(\RR^d)$, so~\eqref{eq:transformer-pde} is generally \emph{not} a Wasserstein gradient flow. This lack of a standard variational structure complicates the analysis of long-time behavior.

\paragraph{Gaussian distributions of tokens.}

To better understand the depth evolution of token distributions, we follow the same approach as in Propositions~\ref{prop:gaussian-diffusion},~\ref{prop:wass-gaussian}, and~\ref{prop:gaussian-mlp}, and specialize to the case where the token distribution is Gaussian. This provides insight into the interplay between the anisotropy of the token distribution and the attention parameters~$\theta$.

\begin{theorem}[Gaussian closure and covariance dynamics~\cite{castin2025unified}]\label{thm:gaussian-closure}
Let the initial distribution be $\alpha_{t=0}=\mathcal N(0,\Sigma_0)$ with covariance $\Sigma_0\succ 0$. Then for as long as a smooth solution to~\eqref{eq:transformer-pde} exists, the distribution remains Gaussian, $\alpha_t=\mathcal N(0,\Sigma_t)$. Moreover, the covariance satisfies the closed ODE
\begin{equation}\label{eq:moment-odes}
    \dot \Sigma_t \;=\; V_t \Sigma_t K_t^\top Q_t \Sigma_t \;+\; \Sigma_t Q_t^\top K_t \Sigma_t V_t^{\top}.
\end{equation}
\end{theorem}

The system~\eqref{eq:moment-odes} is nonlinear and quadratic in $\Sigma_t$. Depending on $(Q,K,V)$ and the initialization $\Sigma_0$, the covariance can blow up in finite time, echoing clustering phenomena observed in related attention-driven dynamics. A detailed analysis—especially for time-dependent parameters $(Q_t,K_t)$ or in settings beyond the potential case—remains an open problem.

\subsection{Training and universality}

The previous section analyzed the depth evolution of an already pretrained transformer. A more difficult question is to understand the \emph{training dynamics}, where one minimizes an empirical loss function. For images, a common choice is a denoising objective~\cite{vincent2011connection,ho2020denoising}, while for text the standard approach is next-token prediction~\cite{radford2019language,brown2020language}. Despite their central importance, the theoretical properties of these training dynamics remain largely open. Empirically, in many large-scale transformer settings, stochastic gradient descent (SGD) alone is less effective than adaptive methods such as Adam~\cite{kingma2015adam}. Adam performs a normalization of the gradient using estimates of its first and second moments, a mechanism that is likely to help stabilize attention layers, which are known to be highly sensitive~\cite{castin2023smooth}.

A related fundamental question is the \emph{expressivity} of transformers. Given depth-dependent parameters $\theta=(\theta(t))_{t\in[0,1]}$, we denote by
\[
   \mathcal{T}_\theta(\alpha_{t=0},x(0)) := x(1)
\]
the transformer mapping defined as the solution at time $t=1$ of the ODE~\eqref{eq:transform-ode-coupled}, starting from the initial distribution $\alpha_{t=0}$ and token $x(0)$. Importantly, $\mathcal{T}_\theta(\alpha,x)$ depends on both the distribution of tokens $\alpha \in \Pp(\RR^d)$ and the individual token position $x\in\RR^d$. Such mappings are often referred to as \emph{in-context mappings}, since the output of a token depends on its surrounding context. Understanding the properties and learnability of in-context mappings has become a central theme in modern machine learning~\cite{garg2022can,xie2021explanation}.
The following universality result, due to Furuya, de Hoop, and Peyré~\cite{furuya2024transformers}, extends the classical universal approximation theorems for MLPs~\cite{cybenko1989approximation,hornik1991approximation} to the transformer setting.

\begin{theorem}[Furuya, de Hoop, Peyré]
Let $\Omega \subset \RR^d$ be compact and let $\mathcal{T}^\star :\Pp(\Omega)\times\Omega \to \RR^d$ be continuous with respect to the product topology combining the Wasserstein metric on $\Pp(\Omega)$ and the Euclidean topology on $\Omega$. Then $\mathcal{T}^\star$ can be uniformly approximated by mappings of the form~$\mathcal{T}_\theta$.
\end{theorem}

Obtaining sharper, quantitative versions of such universality statements is an active area of research; see~\cite{geshkovski2024measure} for recent advances. Another key question is to characterize the structure of optimal depth-dependent parameters $\theta(t)$ for specific prediction tasks, which is essential to understand how in-context learning operates. For example,~\cite{sander2024towards} investigates the ability of transformers to accurately predict elements of recurring sequences in context, by progressively estimating hidden sequence parameters as depth increases.

\section*{Conclusion.}

This survey highlighted two mathematical paradigms that underpin many modern developments in machine learning: diffusion-based models and optimal transport. Both approaches rely on the study of dynamical systems of probability measures, but they emphasize different aspects. Diffusion models, rooted in score estimation and stochastic flows, have become the cornerstone of generative AI thanks to their algorithmic scalability and empirical success. Optimal transport, in contrast, offers a principled geometric framework, enabling variational formulations of nonlinear PDEs and providing tools to analyze convergence, convexity, and stability.
The interplay between these two perspectives reveals a unifying picture: many problems in learning can be recast as the evolution of probability measures, whether describing data distributions in generative modeling, neurons in mean-field limits of neural networks, or tokens in transformers. This dynamical viewpoint sheds light on fundamental questions but also exposes several open challenges.
First, while the theory of Wasserstein flows provides convergence guarantees for certain shallow mean-field models, extending such results to deeper and more expressive architectures remains an outstanding challenge. Second, for diffusion models, much remains to be understood about the structure of the displacement fields they generate: their non-optimality compared to OT, the propagation of discretization errors, and the sample complexity required for accurate score estimation. Finally, transformers introduce an entirely new class of dynamics: their evolution through depth can be modeled by PDEs of interacting particles, but understanding their optimization remains largely open. Key questions include whether one can control such PDEs effectively during training, and how the depth-dependent weights $\theta(t)$ acquire the structure needed for in-context learning.
Bridging these theoretical gaps is not only of mathematical interest but also crucial for the design of future machine learning models. A deeper understanding of these flows may eventually provide both rigorous guarantees and new algorithmic insights, guiding the next generation of architectures that combine the flexibility of diffusion with the structure of optimal transport.


\section*{Acknowledgments.}

A number of the results presented in this survey grew out of stimulating collaborations and exchanges of ideas with colleagues, including Pierre Ablin, Raphaël Barboni, L\'ea Bohbot, Mathieu Blondel, Valérie Castin, Laura Cantini, José Carrillo, Maarten de Hoop, Takashi Furuya, Samuel Hurault, Geert-Jan Huizing, Cyril Letrouit, Thomas Moreau, Jules Samaran, Michael Sander, and François-Xavier Vialard. I am deeply grateful for these joint efforts, which have greatly influenced the perspectives developed here.
I would like to thank Michael Sander for his many comments on this manuscript, which helped improve it.
This work was supported by the European Research Council (ERC project WOLF) and the French government under the management of Agence Nationale de la Recherche as part of the ``France 2030'' program, reference ANR-23-IACL-0008 (PRAIRIE-PSAI).

\bibliographystyle{plain}
\bibliography{biblio}

\end{document}